\documentclass[11pt]{article}
\usepackage[utf8]{inputenc}
\usepackage{amsmath,amssymb,amsthm}
\usepackage{geometry}
\usepackage{hyperref}

\newtheorem{theorem}{Theorem}

\newtheorem{remark}[theorem]{Remark}
\newtheorem{definition}[theorem]{Definition}

\usepackage{mathtools}
\mathtoolsset{showonlyrefs}

\title{A PDE Derivation of the Schr\"odinger--Bass Bridge}

\author{Alexandre Alouadi\thanks{BNP-PAR and Ecole Polytechnique, CMAP}
 \and Pierre Henry-Labord\`ere\thanks{QubeRT}  \and Gr\'egoire Loeper\thanks{BNP-PAR}  \and Othmane Mazhar\thanks{LPSM, Universit\'e Paris Cit\'e}  \and Huy\^en Pham\thanks{Ecole Polytechnique, CMAP}  \and Nizar Touzi\thanks{NYU Tandon school of Engineering}}

\date{\today}

\begin{document}

\maketitle

\begin{abstract}
This short paper announces the main results of \cite{SBB2026}, where the Schr\"odinger--Bass Bridge (SBB) problem is introduced and studied in full generality. Here we provide a direct PDE derivation of the SBB system in dimension one, showing how the optimal coupling problem that interpolates between the classical Schr\"odinger bridge and the Bass martingale transport can be solved explicitly via Legendre transforms and the heat equation. A key insight is that the optimal SBB process is a \emph{Stretched Schr\"odinger Bridge}: the composition of a monotone transport map with a Schr\"odinger bridge. This extends the stretched Brownian motion representation of Bass martingales to the semimartingale setting and provides a unified framework that recovers both the Sinkhorn algorithm (in the limit $\beta \to \infty$) and the Bass construction (as $\beta \to 0$). We refer to \cite{SBB2026} for complete proofs, the multidimensional setting, strong duality, dual attainment, and further developments.
\end{abstract}

\section{Introduction}

The Schr\"odinger bridge problem and martingale optimal transport provide two fundamental approaches to constructing optimal couplings of probability measures via continuous-time stochastic processes.

\paragraph{The Schr\"odinger bridge.} Given a reference Brownian motion and prescribed initial and terminal distributions $\mu_0$ and $\mu_T$, the Schr\"odinger problem seeks the path measure that minimizes relative entropy with respect to the Brownian reference while satisfying the marginal constraints. This problem has deep connections to large deviation theory and has recently found applications in machine learning and generative modeling; see \cite{Leonard2014,Chen2021} and the references therein.

\paragraph{The Bass martingale.} In contrast, the Bass problem (or stretched Brownian motion) prescribes the marginal laws of a continuous martingale and minimizes a pathwise cost measuring deviation from Brownian behavior. The solution can be represented as a monotone transport of Brownian motion, leading to elegant probabilistic constructions; see \cite{BachhoffBeigblock2020,BachhoffSchachermayer2025,ConzeHL2021,JosephLoeperObloj2024}.

\paragraph{The Schr\"odinger--Bass Bridge.} The work \cite{SBB2026} introduces a one-parameter family of optimal semimartingale couplings that interpolates between these two classical problems. For a parameter $\beta > 0$ and probability measures $\mu_0, \mu_T \in \mathcal{P}_2(\mathbb{R}^d)$, one considers continuous semimartingales $X = (X_t)_{0 \leq t \leq T}$ with $X_0 \sim \mu_0$, $X_T \sim \mu_T$, and minimizes a functional that penalizes both the drift and the volatility of $X$:
\begin{equation}\label{eq:SBB-cost}
\text{SBB}(\mu_0, \mu_T) := \inf_{(X_t)} \mathbb{E}\left[\frac{1}{2}\int_0^T \left(|\alpha_t|^2 + \beta|\sigma_t - \text{Id}|^2\right) dt\right],
\end{equation}
where $\alpha_t$ and $\sigma_t$ denote the drift and volatility of the semimartingale $X$.

The key insight is that:
\begin{itemize}
\item As $\beta \to \infty$, the volatility is constrained to be close to the identity, and one recovers the classical Schr\"odinger bridge.
\item As $\beta \to 0$, the volatility penalization vanishes; the optimal process becomes a martingale, and one recovers the Bass problem.
\end{itemize}

\paragraph{Contribution.} In this note, we provide a direct PDE derivation of the main structural results of \cite{SBB2026} in dimension one. By working through the Hamilton--Jacobi--Bellman equation and exploiting Legendre duality, we show that the optimal time-marginals admit a representation involving:
\begin{enumerate}
\item A function $h$ solving the backward heat equation,
\item A measure $\nu$ evolving by the forward heat equation,
\item A pair of inverse transport maps $\mathcal{Y}$ and $\mathcal{X} = \mathcal{Y}^{-1}$ given by gradients of convex functions.
\end{enumerate}
This yields the \emph{SBB system}, which generalizes both the classical Sinkhorn iterations and the martingale Sinkhorn algorithm of \cite{JosephLoeperObloj2024}.

A key conceptual insight is that the process $Y_t := \mathcal{Y}(t, X_t)$ is itself a \emph{Schr\"odinger bridge}, and the optimal SBB process $X_t = \mathcal{X}(t, Y_t)$ is obtained by applying a monotone transport map to it. We call this a \textbf{Stretched Schr\"odinger Bridge}, by analogy with the stretched Brownian motion representation of Bass martingales. This provides a unified perspective:
\begin{center}
\emph{Bass martingale = stretched Brownian motion}\\
\emph{SBB solution = stretched Schr\"odinger bridge}. 
\end{center}

\section{Setup and Dual Problem}

We work in dimension $d = 1$ for notational simplicity; the derivation extends verbatim to any dimension $d \geq 1$. Consider the SDE
\begin{equation}\label{eq:SDE}
dX_t = \alpha_t \, dt + \sigma_t \, dW_t, \quad X_0 \sim \mu_0,
\end{equation}
where $(\alpha_t, \sigma_t)_{0 \leq t \leq T}$ are adapted processes. The \textbf{primal problem} is
\begin{equation}\label{eq:primal}\tag{P}
\text{SBB}(\mu_0, \mu_T) := \inf_{(\alpha, \sigma) \text{ adapted}} \mathbb{E}\left[\frac{1}{2} \int_0^T \left(\alpha_t^2 + \beta(\sigma_t - 1)^2\right) dt\right]
\end{equation}
subject to the marginal constraints $X_0 \sim \mu_0$ and $X_T \sim \mu_T$.

\subsection{The dual problem}

By standard arguments in stochastic control (see \cite{TanTouzi2013}), the \textbf{dual problem} takes the form
\begin{equation}\label{eq:dual}\tag{D}
\mathcal{V}(\mu_0, \mu_T) := \sup_v \left\{ \mathbb{E}_{\mu_T}[v(T, X_T)] - \mathbb{E}_{\mu_0}[v(0, X_0)] \right\}
\end{equation}
over functions $v$ satisfying the Hamilton--Jacobi--Bellman (HJB) equation
\begin{equation}\label{eq:HJB-u}
\partial_t v + \frac{1}{2}|\partial_x v|^2 + \frac{1}{2}\frac{\partial_{xx} v}{1 - \partial_{xx}v/\beta} = 0,
\end{equation}
with the constraint that $\partial_{xx} v < \beta$ (equivalently, the function $u := x^2/2 - v/\beta$ is strictly convex).

\begin{remark}[Strong duality and attainment]
The proof of strong duality $\text{SBB}(\mu_0, \mu_T) = \mathcal{V}(\mu_0, \mu_T)$ and the attainment of the infimum in the primal problem \eqref{eq:primal} are studied in detail in \cite{SBB2026}. This is a delicate matter: standard results from semimartingale optimal transport (e.g., \cite{TanTouzi2013}) do not immediately apply because the running cost $(\sigma - 1)^2$ lacks coercivity in $\sigma$. The analysis in \cite{SBB2026} overcomes this difficulty through careful a priori estimates and a reduction to a static dual formulation involving $\beta$-convex potentials.
\end{remark}

\section{PDE Derivation of the SBB System}

\subsection{Change of variables}

Define the convex function
\[
u := \frac{x^2}{2} - \frac{v}{\beta}.
\]
Substituting into \eqref{eq:HJB-u}, the function $u$ satisfies
\begin{equation}\label{eq:HJB-v}
\partial_t u - \frac{\beta}{2}|\partial_x u - x|^2 + \frac{1}{2}\left(1 - \frac{1}{\partial_{xx}u}\right) = 0.
\end{equation}

\subsection{Legendre transform}

Let $u^*$ denote the Legendre transform of $u$ in the spatial variable:
\[
u^*(t, y) := \sup_{x \in \mathbb{R}} \left\{ xy - u(t, x) \right\}.
\]
Since $u$ is convex in $x$, a standard computation shows that $u^*$ satisfies
\begin{equation}\label{eq:HJB-vstar}
\partial_t u^* + \frac{\beta}{2}|\partial_y u^* - y|^2 - \frac{1}{2}(1 - \partial_{yy}u^*) = 0.
\end{equation}

\subsection{Reduction to the heat equation}

Define
\[
w := u^* - \frac{y^2}{2}.
\]
Substituting into \eqref{eq:HJB-vstar}, the function $w$ satisfies
\begin{equation}\label{eq:HJB-w}
\partial_t w + \frac{\beta}{2}|\partial_y w|^2 + \frac{1}{2}\partial_{yy}v^* = 0.
\end{equation}
Remarkably, if we set
\begin{equation}\label{eq:h-def}
h := \exp(\beta w), \quad \text{equivalently} \quad w = \frac{1}{\beta}\log h,
\end{equation}
then $h$ solves the \textbf{heat equation}:
\begin{equation}\label{eq:heat}
\partial_t h + \frac{1}{2}\partial_{yy}h = 0.
\end{equation}
This is the key observation: the nonlinear HJB structure reduces to the linear heat equation after the appropriate sequence of transforms.

\section{The Transport Maps $\mathcal{Y}$ and $\mathcal{X}$}

A key structural feature of the SBB problem is that the optimal solution can be encoded via a pair of inverse transport maps $\mathcal{Y}$ and $\mathcal{X}$.

\subsection{Definition of the maps}

Since $u$ is strictly convex in $x$, the function $y \mapsto \log h(t,y) + \frac{\beta}{2}|y|^2$ is strictly convex (recall $h = e^{\beta w}$ and $w = u^* - y^2/2$). For each $(t,x)$, define
\begin{equation}\label{eq:Y-map}
\boxed{\mathcal{Y}(t, x) := \arg\min_{y \in \mathbb{R}} \left\{ \log h(t, y) + \frac{\beta}{2}|x - y|^2 \right\}.}
\end{equation}
The first-order condition for this minimization gives
\[
\partial_y \log h(t, \mathcal{Y}(t,x)) + \beta(\mathcal{Y}(t,x) - x) = 0,
\]
which can be rewritten as
\begin{equation}\label{eq:Y-implicit}
\mathcal{Y}(t, x) = x - \frac{1}{\beta}\partial_y \log h(t, \mathcal{Y}(t, x)).
\end{equation}

Define the map $\mathcal{X}$ by 
\begin{equation}\label{eq:X-map}
\boxed{\mathcal{X}(t, y) := y + \frac{1}{\beta}\partial_y \log h(t, y).}
\end{equation}
One readily verifies that $\mathcal{X}(t, \mathcal{Y}(t, x)) = x$ and $\mathcal{Y}(t, \mathcal{X}(t, y)) = y$, i.e.,  $\mathcal{X} = \mathcal{Y}^{-1}$ is the inverse map (w.r.t. the spatial variable)  of $\mathcal{Y}$. 

\begin{remark}[Gradient structure]
The map $\mathcal{X}$ is the gradient of a convex function:
\[
\mathcal{X}(t, y) = \partial_y \mathcal{G}(t, y), \quad \text{where} \quad \mathcal{G}(t, y) := \frac{y^2}{2} + \frac{1}{\beta}\log h(t, y) = u^*(t, y).
\]
\end{remark}

\subsection{Connection with the dual variables}

The maps $\mathcal{Y}$ and $\mathcal{X}$ are directly related to the gradient of the dual value function $v$:
\begin{equation}\label{eq:Y-via-u}
\mathcal{Y}(t, x) = x - \frac{1}{\beta}\partial_x v(t, x).
\end{equation}
Indeed, from the Legendre transform relations $v = \beta(x^2/2 - u)$ and $u^* = y^2/2 + w$, one can verify that the minimizer in \eqref{eq:Y-map} satisfies \eqref{eq:Y-via-u}.

\section{The Optimal Process and the SBB System}

\subsection{Dynamics of the optimal process}

Let $X_t$ be the process induced by the optimal control from the solution $u$ to the HJB equation, i.e.,
\begin{equation}\label{eq:optimal-SDE}
dX_t = \partial_x v(t,X_t) \, dt + \frac{1}{1 - \partial_{xx}v(t,X_t)/\beta} \, dW_t.
\end{equation}
A standard verification argument shows that $\partial_x v(t, X_t)$ is a martingale.

Using the transport maps, the optimal drift and volatility can be expressed as:
\begin{align}
\alpha(t, x) &= \partial_x v(t, x) = \beta\big(x - \mathcal{Y}(t, x)\big) = \partial_y \log h\big(t, \mathcal{Y}(t, x)\big), \label{eq:alpha-via-Y}\\
\sigma(t, x) &= \frac{1}{1 - \partial_{xx}v(t,x)/\beta} = 1 + \frac{1}{\beta}\partial_{yy}\log h\big(t, \mathcal{Y}(t, x)\big). \label{eq:sigma-via-Y}
\end{align}

\subsection{The transformed process}

Define the process
\begin{equation}\label{eq:Y-def}
Y_t := \mathcal{Y}(t, X_t) = X_t - \frac{1}{\beta}\partial_x v(t, X_t).
\end{equation}
By It\^o's formula, one can show that
\[
dY_t = \partial_x v(t, X_t) \, dt + dW_t.
\]
Moreover, by the chain of Legendre transforms, we have
\begin{align*}
\partial_x v(t, x) &= \beta(x - \partial_x u) \; = \;  \beta(\partial_y u^* - y) \\
& = \;  \beta \partial_y w \; = \;  \partial_y \log h.
\end{align*}
Therefore,
\begin{equation}\label{eq:Y-SDE}
dY_t = \partial_y \log h(t, Y_t) \, dt + dW_t.
\end{equation}

\subsection{Inverse map and stretched Brownian motion}

The inverse relationship, using the map $\mathcal{X}$, is
\begin{equation}\label{eq:X-from-Y}
X_t = \mathcal{X}(t, Y_t) = Y_t + \frac{1}{\beta}\partial_y \log h(t, Y_t).
\end{equation}
This shows that the optimal process $X$ is obtained by applying the time-dependent gradient map $\mathcal{X}$ to the process $Y$.

\section{Stretched Schr\"odinger Bridge}

We now make a key observation: the process $Y_t$ is itself a \emph{Schr\"odinger bridge}, and the optimal SBB process $X_t$ is obtained by applying a monotone transport to it. We first recall a standard characterization of Schr\"odinger bridges (see \cite{Leonard2014, Chen2021} for detailed treatments).

\begin{definition}[Schr\"odinger Bridge]\label{def:SB}
Let $h : [0,T] \times \mathbb{R} \to (0,\infty)$ be a positive solution of the backward heat equation $\partial_t h + \frac{1}{2}\partial_{yy}h = 0$. A process $(Y_t)_{0 \leq t \leq T}$ is called a \textbf{Schr\"odinger bridge} (with potential $h$) if it satisfies
\begin{equation}\label{eq:SB-SDE}
dY_t = \partial_y \log h(t, Y_t) \, dt + dW_t.
\end{equation}
Equivalently, there exists a probability measure $\mathbb{Q}$ under which $Y$ is a standard Brownian motion, and $\mathbb{P}$ is related to $\mathbb{Q}$ by
\[
\frac{d\mathbb{P}}{d\mathbb{Q}} = h(T, Y_T).
\]
\end{definition}

The drift $\partial_y \log h$ in \eqref{eq:SB-SDE} is the classical \emph{Nelson drift} or \emph{osmotic velocity} that characterizes Schr\"odinger bridges. The function $h$ is space-time harmonic for the backward heat operator.

\begin{definition}[Stretched Schr\"odinger Bridge]\label{def:SSB}
Let $(Y_t)_{0 \leq t \leq T}$ be a Schr\"odinger bridge with potential $h$, and let $\mathcal{X} : [0,T] \times \mathbb{R} \to \mathbb{R}$ be a family of monotone increasing maps such that $\mathcal{X}(t, \cdot) = \partial_y \mathcal{G}(t, \cdot)$ for some convex function $\mathcal{G}(t, \cdot)$. The process
\begin{equation}\label{eq:SSB-def}
X_t := \mathcal{X}(t, Y_t)
\end{equation}
is called a \textbf{Stretched Schr\"odinger Bridge}.
\end{definition}

\begin{remark}[The stretching potential $\mathcal{G}$]\label{rem:stretching-potential}
In Definition \ref{def:SSB}, the stretching map $\mathcal{X}(t, \cdot) = \partial_y \mathcal{G}(t, \cdot)$ is the gradient of a convex potential $\mathcal{G}$. A key feature of the Stretched Schr\"odinger Bridge is that \textbf{the stretching potential $\mathcal{G}$ is directly determined by the Schr\"odinger potential $h$} of the underlying Schr\"odinger bridge:
\begin{equation}\label{eq:G-from-h}
\mathcal{G}(t, y) = \frac{|y|^2}{2} + \frac{1}{\beta}\log h(t, y).
\end{equation}
The first term $\frac{|y|^2}{2}$ is the identity potential, and the second term $\frac{1}{\beta}\log h$ is a perturbation controlled by:
\begin{itemize}
\item the Schr\"odinger potential $h$, which encodes the drift of the base Schr\"odinger bridge $Y$;
\item the parameter $\beta$, which controls the strength of the perturbation.
\end{itemize}
Thus, the stretching map
\[
\mathcal{X}(t, y) = \partial_y \mathcal{G}(t, y) = y + \frac{1}{\beta}\partial_y \log h(t, y) = y + \frac{1}{\beta}\frac{\partial_y h(t,y)}{h(t,y)}
\]
combines the identity with the normalized gradient of the Schr\"odinger potential $h$.
\end{remark}

\begin{remark}[Hierarchy of processes]
The Stretched Schr\"odinger Bridge reveals a natural hierarchy of optimal transport constructions:
\begin{center}
\begin{tabular}{cccc}
\textbf{Base process} & \textbf{Potential $h$} & \textbf{Stretching potential $\mathcal{G}$} & \textbf{Result} \\
\hline
Brownian motion $W$ & $h \equiv 1$ & $\mathcal{G}(y) = \frac{|y|^2}{2} + \frac{1}{\beta}\psi(y)$ & Bass martingale \\[2pt]
Schr\"odinger bridge $Y$ & $h$ (heat kernel) & $\mathcal{G}(y) = \frac{|y|^2}{2}$ & Schr\"odinger bridge \\[2pt]
Schr\"odinger bridge $Y$ & $h$ (heat kernel) & $\mathcal{G}(y) = \frac{|y|^2}{2} + \frac{1}{\beta}\log h$ & Stretched SB
\end{tabular}
\end{center}
In the Bass martingale case \cite{BachhoffSchachermayer2025}, the base process is Brownian ($h \equiv 1$) and the stretching potential $\mathcal{G}$ is determined by a separate convex function $\psi$. In the Stretched Schr\"odinger Bridge, both the base Schr\"odinger bridge and the stretching map are coupled through the \emph{same} potential $h$.
\end{remark}

\begin{remark}[Connection to optimal transport calibration]\label{rem:calibration}
The Stretched Schr\"odinger Bridge framework has direct applications to \emph{model calibration} in mathematical finance. In the calibration problem, one seeks to construct a diffusion process $X$ whose marginal distributions match observed option prices (via the Breeden--Litzenberger formula), while staying close to some reference model.

The Bass martingale approach to calibration was developed in \cite{ConzeHL2021}, where the stretched Brownian motion construction provides a local volatility model consistent with given marginals. The path-dependent optimal transport perspective on calibration was studied in \cite{GuoLoeper2021}, connecting optimal transport maps to model calibration on exotic derivatives. See also \cite{AcciaioMariniPammer2023} for calibration of the Bass local volatility model.

The SBB framework extends these approaches by:
\begin{enumerate}
\item allowing a non-trivial reference model (the Schr\"odinger bridge component);
\item providing a principled way to balance fit to marginals (via the stretching map $\mathcal{X}$) against proximity to the reference (via the parameter $\beta$);
\item yielding explicit formulas for the optimal drift and volatility in terms of the Schr\"odinger potential $h$.
\end{enumerate}
\end{remark}

\begin{theorem}[SBB solution as Stretched Schr\"odinger Bridge]\label{thm:SSB}
The optimal process $X_t$ solving the SBB problem \eqref{eq:primal} is a Stretched Schr\"odinger Bridge. Specifically:
\begin{enumerate}
\item The process $Y_t = \mathcal{Y}(t, X_t)$ is a Schr\"odinger bridge with potential $h = e^{\beta w}$, where $h$ solves the backward heat equation $\partial_t h + \frac{1}{2}\partial_{yy}h = 0$ with terminal condition  determined by the marginal constraints.
\item The optimal process is recovered as
\[
X_t = \mathcal{X}(t, Y_t) = Y_t + \frac{1}{\beta}\partial_y \log h(t, Y_t),
\]
where $\mathcal{X}(t, \cdot) = \partial_y \mathcal{G}(t, \cdot)$ with the convex potential
\[
\mathcal{G}(t, y) = \frac{y^2}{2} + \frac{1}{\beta}\log h(t, y) = u^*(t, y).
\]
\item There exists a probability measure $\mathbb{Q} \sim \mathbb{P}$ with $\frac{d\mathbb{P}}{d\mathbb{Q}} = h(T, Y_T)$ under which $Y$ is a standard Brownian motion.
\end{enumerate}
\end{theorem}

\begin{remark}[Interpolation via the stretching map]
The parameter $\beta$ controls the ``amount of stretching'':
\begin{itemize}
\item As $\beta \to \infty$: $\mathcal{X}(t, y) \to y$, so $X_t \to Y_t$. The stretched Schr\"odinger bridge reduces to an ordinary Schr\"odinger bridge.
\item As $\beta \to 0$: $h \to \text{const}$, so $Y_t$ becomes a Brownian motion. The stretched Schr\"odinger bridge reduces to a stretched Brownian motion (Bass martingale).
\end{itemize}
Thus, the Stretched Schr\"odinger Bridge provides a unified framework that interpolates between the Schr\"odinger bridge ($\beta = \infty$) and the Bass martingale ($\beta = 0$).
\end{remark}

\section{The SBB System}

Combining the above, we obtain the \textbf{Schr\"odinger--Bass Bridge system}. The key observation is that the time-marginals $\mu_t$ of the optimal process $X$ can be expressed using the transport map $\mathcal{X}$:
\begin{equation}\label{eq:SBB-system}
\boxed{
\begin{aligned}
&X_t \sim \mu_t, \\
&Y_t \sim h(t,\cdot)  \nu_t, \\
&\mu_t = \mathcal{X}(t, \cdot)_\# (h(t,\cdot) \nu_t), \\
&\partial_t h + \frac{1}{2}\partial_{xx}h = 0, \\
&\partial_t \nu - \frac{1}{2}\partial_{xx}\nu = 0,
\end{aligned}
}
\end{equation}
where $\mathcal{X}(t, y) = y + \frac{1}{\beta}\partial_y \log h(t,y)$ is the transport map  in \eqref{eq:X-map}, and $\nu$ is the density of the Brownian motion starting from some initial distribution $\nu_0$, and satisfying the Fokker-Planck equation. 

Equivalently, using the inverse map $\mathcal{Y} = \mathcal{X}^{-1}$, the system can be written as
\begin{equation}\label{eq:SBB-system-Y}
\begin{cases}
\mathcal{Y}(T, \cdot)_\# \mu_T = h(T,\cdot)  \nu_T, \\
\mathcal{Y}(0, \cdot)_\# \mu_0 = h(0,\cdot)  \nu_0.
\end{cases}
\end{equation}
This is the \emph{Schr\"odinger--Bass system} of \cite{SBB2026}, which characterizes the optimal solution through the boundary conditions at times $0$ and $T$.

\begin{remark}[Interpolation between Sinkhorn and Bass]
The SBB system interpolates between two classical algorithms:

\textbf{Sinkhorn ($\beta \to \infty$):} As $\beta \to \infty$, the maps $\mathcal{X}$ and $\mathcal{Y}$ both converge to the identity:
\[
\mathcal{X}(t, y) = y + \frac{1}{\beta}\partial_y \log h(t,y) \to y, \quad \mathcal{Y}(t, x) \to x.
\]
We recover the classical Schr\"odinger system:
\[
X_t \sim \mu_t, \quad Y_t \sim h(t,\cdot)  \nu_t, \quad \mu_t = h(t,\cdot)  \nu_t,
\]
with $X_t = Y_t$.

\textbf{Bass ($\beta \to 0$):} As $\beta \to 0$, the function $h$ tends to a constant (say $h \equiv 1$), and the map $\mathcal{X}$ becomes purely determined by the convex potential $u^*$:
\[
\mathcal{X}(t, y) = \partial_y u^*(t, y).
\]
We recover the Bass martingale construction:
\[
X_t \sim \mu_t, \quad Y_t \sim \nu_t, \quad \mu_t = (\partial_y u^*(t,\cdot))_\# \nu_t,
\]
where $X = \partial_y u^*(Y)$ for $Y$ Brownian.
\end{remark}

\section{Sinkhorn-Type Iterations}

To solve the SBB system numerically, one needs to enforce the marginal constraints at times $0$ and $T$. Using the map $\mathcal{Y}$, the terminal constraint \eqref{eq:SBB-system-Y} becomes
\begin{equation}\label{eq:terminal-constraint}
\mathcal{Y}(T, \cdot)_\# \mu_T = h(T,\cdot)  \nu_T = e^{\beta w(T,\cdot)} \nu_T,
\end{equation}
recalling that $h = e^{\beta w}$ where $w = u^* - y^2/2$.

Equivalently, using the inverse map $\mathcal{X}$, we need
\[
\mu_T = \mathcal{X}(T, \cdot)_\# (e^{\beta w(T,\cdot)} \nu_T).
\]
With $u^* = y^2/2 + w$ at time $T$, this becomes
\[
\mu_T = (\partial_y u^*(T,\cdot))_\# \left(\exp\left(\beta\left(u^*(T,\cdot) - \frac{y^2}{2}\right)\right) \nu_T\right).
\]

This leads to the Monge--Amp\`ere equation
\begin{equation}\label{eq:MA}
\det D^2 u^*(T,\cdot) = \frac{\nu}{\mu(\partial_y u^*(T,\cdot))} \exp\left(\beta\left(u^*(T,\cdot) - \frac{y^2}{2}\right)\right).
\end{equation}

\begin{remark}[Structure of the Monge--Amp\`ere equation]\label{rem:MA}
The equation \eqref{eq:MA} exhibits a remarkable structure that blends two classical Monge--Amp\`ere theories:

\textbf{(i) Optimal transport.} The standard Monge--Amp\`ere equation arising in optimal transport with quadratic cost takes the form
\[
\det D^2 \upsilon = \frac{f}{g(\nabla \upsilon)},
\]
where $f$ and $g$ are the source and target densities, and $\nabla \upsilon$ is the optimal transport map. This equation, studied extensively by Brenier \cite{Brenier1991} and Caffarelli \cite{Caffarelli1992}, characterizes the convex potential whose gradient pushes forward $f$ to $g$. Regularity theory for such equations, including the Ma--Trudinger--Wang (MTW) curvature conditions, was developed in \cite{MaTrudingerWang2005,Loeper2009}.

\textbf{(ii) K\"ahler geometry.} In K\"ahler geometry and the study of canonical metrics, one encounters Monge--Amp\`ere equations of the form
\[
\det D^2 \upsilon = e^{f + \lambda \upsilon},
\]
where the right-hand side depends exponentially on the potential $u$ itself. The case $\lambda > 0$ arises in K\"ahler--Einstein metrics with positive Ricci curvature (Aubin \cite{Aubin1976}, Yau \cite{Yau1978}), while the affine sphere equation $\det D^2 \upsilon = e^{\pm \upsilon}$ was studied by Calabi and Cheng--Yau \cite{ChengYau1986}.

\textbf{(iii) The SBB equation as a blend.} The Monge--Amp\`ere equation \eqref{eq:MA} combines both structures:
\[
\det D^2 u^*(T,\cdot) = \underbrace{\frac{\nu}{\mu(\partial_y u^*(T,\cdot))}}_{\text{OT Jacobian}} \times \underbrace{\exp\left(\beta\left(u^*(T,\cdot) - \frac{|y|^2}{2}\right)\right)}_{\text{Entropic/K\"ahler weight}}.
\]
The first factor is the classical optimal transport Jacobian, enforcing the pushforward constraint $(\partial_y u^*(T,\cdot))_\# \nu = \mu$. The second factor is an exponential weight depending on the potential $u^*(T,\cdot)$ itself, akin to the K\"ahler--Einstein and affine sphere equations. The parameter $\beta$ interpolates between these two regimes:
\begin{itemize}
\item As $\beta \to 0$: the exponential weight becomes constant, and we recover the pure optimal transport Monge--Amp\`ere equation.
\item As $\beta \to \infty$: the exponential weight dominates, enforcing $u^*(T,\cdot) \approx |y|^2/2$ (i.e., $\partial_y u^*(T,\cdot) \approx \text{Id}$), and the equation degenerates to the Schr\"odinger system.
\end{itemize}
This blend suggests that regularity and existence theory for \eqref{eq:MA} may require techniques from both optimal transport (MTW conditions, Caffarelli's regularity) and K\"ahler geometry (a priori estimates, continuity methods). We leave the detailed analysis of this equation to future work.
\end{remark}

An iterative scheme alternating between:
\begin{enumerate}
\item Updating the spatial map $\mathcal{Y}(0, x) = \arg\min_y \left[\log h(0, y) + \frac{\beta|x-y|^2}{2}\right]$,
\item Updating $w(T,\cdot)$ via the log-density ratio $w^{k+1}(T,\cdot) = \frac{1}{\beta}\log \frac{d(\mathcal{Y}^k(T,\cdot)_\# \mu_T)}{d\nu_T^k}$,
\end{enumerate}
provides a Sinkhorn-like algorithm for computing the SBB solution. The maps $\mathcal{Y}$ and $\mathcal{X}$ are then recovered from the converged potential via \eqref{eq:Y-map} and \eqref{eq:X-map}. 


\end{document}